\documentclass{article}
\usepackage[paperwidth=7in, paperheight=10in, margin=.875in]{geometry}
\usepackage{amsfonts,amssymb}
\usepackage{amsmath}
\usepackage{graphicx}
\usepackage{cite}
\usepackage{enumerate}
\usepackage{color}
\numberwithin{equation}{section}
\newtheorem{theorem}{Theorem}%[section]

\newtheorem{lemma}{Lemma}%[section]
\newtheorem{proposition}{Proposition}%[section]
\newtheorem{example}{Example}%[section]
\newtheorem{remark}{Remark}%[section]

\def\begproof{\noindent{\bf Proof: }}
\def\endproof{\par\rightline{\vrule height5pt width5pt depth0pt}\medskip}

\usepackage[colorinlistoftodos]{todonotes}
\usepackage[cal=boondox]{mathalfa}
\usepackage{mathtools}
\mathtoolsset{showonlyrefs}

\newcommand{\R}{\mathbb{R}}
\newcommand{\N}{\mathbb{N}}
\newcommand{\eps}{\epsilon}

\newcommand{\argmin}[0]{\operatorname{argmin}}
\newcommand{\sign}[0]{\operatorname{sign}}

\def\({\begin{eqnarray}}
\def\){\end{eqnarray}}
\def\[{\begin{eqnarray*}}
\def\]{\end{eqnarray*}}

\def\d{\,\mathrm{d}}

   \allowdisplaybreaks
\begin{document}

\centerline{{\Large\textbf{A mesoscopic model of biological transportation networks}}}
\vskip 7mm

%%%%%%%%%%%%%%%%%%%%%%%%%%%%%%%%%
%% Authors
%%%%%%%%%%%%%%%%%%%%%%%%%%%%%%%%%
\centerline{
	{\large Martin Burger}\footnote{Institut f\"ur Numerische und Angewandte Mathematik and Cells in Motion Cluster of Excellence, Westf\"alische Wilhelms Universit\"at M\"unster, Einsteinstrasse 62, 48149 M\"unster, D; 
		{\it martin.burger@wwu.de}},
	{\large Jan Haskovec}\footnote{Mathematical and Computer Sciences and Engineering Division,
		King Abdullah University of Science and Technology,
		Thuwal 23955-6900, Kingdom of Saudi Arabia; 
		{\it jan.haskovec@kaust.edu.sa}},
}
\centerline{
	{\large Peter Markowich}\footnote{Mathematical and Computer Sciences and Engineering Division,
		King Abdullah University of Science and Technology,
		Thuwal 23955-6900, Kingdom of Saudi Arabia, and University of Vienna, Faculty of Mathematics, Oskar-Morgenstern-Platz 1, 1090 Wien, AT;
		{\it peter.markowich@kaust.edu.sa, peter.markowich@univie.ac.at}},
	{\large Helene Ranetbauer}\footnote{University of Vienna, Faculty of Mathematics, Oskar-Morgenstern-Platz 1, 1090 Wien, AT;
		{\it helene.ranetbauer@univie.ac.at}}
	}
\vskip 10mm

\noindent{\bf Abstract.}
          We introduce a mesoscopic model for natural network formation processes, acting as a bridge between the discrete and continuous network approach proposed in \cite{hu2013adaptation}.
	  The models are based on a common approach where the dynamics of the conductance network is subject to pressure force effects.
	  We first study topological properties of the discrete model and we prove that if the metabolic energy consumption term is concave with respect
	  to the conductivities, the optimal network structure is a tree (i.e., no loops are present).
	  We then analyze various aspects of the mesoscopic modeling approach, in particular its relation to the discrete model and its stationary solutions,
	  including discrete network solutions. Moreover, we present an alternative formulation of the mesoscopic model that avoids the explicit presence of the pressure
	  in the energy functional.
	  \vskip 7mm

%\allowdisplaybreaks

\noindent{\bf AMSC:} 35B36; 92C42; 35K55; 49J20
\vspace{2mm}


%%%%%%%%%%%%%%%%%%%%%%%%%%%%%%%%%%%%%%%%%%%%%%%%%%%%%%%%%%%%%%%%%%%%%%%%%%%%%%%%%%%%%%%%%%%%
\section{Introduction}\label{intro}
%%%%%%%%%%%%%%%%%%%%%%%%%%%%%%%%%%%%%%%%%%%%%%%%%%%%%%%%%%%%%%%%%%%%%%%%%%%%%%%%%%%%%%%%%%%%
Transportation networks play a fundamental role in biological applications such as leaf venation in plants \cite{malinowski2013understanding},
vascular pattern formation \cite{sedmera2011function}, mammalian circulatory systems that convey nutrients to the body through blood circulation, or
neural networks that transport electric charge \cite{eichmann2005guidance, michel1995morphogenesis}.
They have been widely investigated by the scientific community and different tools for describing their development, function and adaptation have been proposed in the literature.

Our work is based on the discrete modeling approach introduced in \cite{hu2013optimization} and  \cite{hu2013adaptation}.
The authors proposed a purely local dynamic adaptation model based on mechanical laws, consisting of a system of ordinary differential equations (ODE) on graph edges,
coupled to a linear system of equations (Kirchhoff law). This system is obtained as the gradient flow of an energy functional consisting of a kinetic energy term (pumping power)
and a metabolic cost term, written in terms of the edge conductivities.
The first contribution of this paper is a proof that the global energy minimizer for the energy with concave metabolic cost term does not contain any loops,
i.e., it is a tree (in the graph-theoretical sense). %We also provide an example showing that the minimizer with convex metabolic cost may or may not contain loops.
%In particular, the model responds only to local information and fluctuations in flow distributions can be naturally incorporated.

In \cite{hu2013optimization} a related PDE-based continuum model was proposed,
which was subsequently studied in the series of papers \cite{haskovec2015mathematical, haskovec2016notes, albi2016biological, albi2017continuum}.
The continuum model consists of a parabolic reaction-diffusion equation for the conductivity vector, constrained by a Poisson equation for the pressure.
Again, the system possesses a constrained gradient flow structure with respect to the continuum version of the energy functional.
In this paper we study the mesoscopic model briefly introduced in \cite{albi2017continuum} as a bridge between the microscopic and macroscopic descriptions.
The model has a formal Wasserstein-type gradient flow structure, constrained again by a Poisson equation.
We analyze its properties and discuss several special stationary solutions including discrete network solutions.
Since the discrete network solution is an inherently one-dimensional structure, its embedding into the multiple-dimensional space
requires a suitable interpretation of the Poisson equation with measure valued permeability tensor.
We formally introduce such an interpretation, showing that the Poisson equation then reduces to a coupled system
of elliptic equations posed on the edges of the network.
Moreover, we show that in the particular regime where edges of the network are aligned with pressure gradients,
the model reduces to the system studied by Putti et al. \cite{facca2018towards}.
This system is a continuous version of the discrete model \cite{Tero} used to simulate the
ability of the slime mold (Physarum Polycephalum) to find the shortest path connecting two food sources in a maze.

To overcome possible solvability issues in the Poisson equation due to the degeneracy of the permeability tensor,
we shall introduce the pressureless formulation of the model, where the Poisson equation is replaced by a linear constraint on the flux variable.
We also reformulate the energy functional, using the formula by Benamou and Brenier (see, e.g., see \cite{santambrogio2015optimal}),
such that its domain of definition can be extended to the space of measures.
We then formulate a minimizing movement scheme with respect to the Fisher-Rao distance, which is expected to provide
solutions of the transient system. However, a rigorous passage to the limit remains an open problem
due to the lack of regularity estimates.

This paper is organized as follows. In Section \ref{sec:discrete} we present the discrete model introduced in \cite{hu2013adaptation} and study its topological properties,
in particular, the presence or absence of loops in the global minimizer of the energy functional.
In Section \ref{sec:mesoscopic} we introduce a mesoscopic modeling approach representing the evolution of a probability measure
to have an edge in a certain direction of some conductivity at a specific point in space and time.
We analyze various aspects including its relation to the discrete model and its stationary solutions.
In Section \ref{sec:monokinetic} we deduce an evolution equation for the contuctivities from the mesoscopic approach via a monokinetic ansatz.
Finally, in Section \ref{sec:pressureless} we present an alternative formulation of the model %introduced in Section \ref{sec:monokinetic}
avoiding the explicit presence of the pressure and formulate the minimizing movement scheme.

%%%%%%%%%%%%%%%%%%%%%%%%%%%%%%%%%%%%%%%%%%%%%%%%%%%%%%%%%%%%%%%%%%%%%%%%%%%%%%%%%%%%%%%%%%%%
\section{The discrete model}\label{sec:discrete}
%%%%%%%%%%%%%%%%%%%%%%%%%%%%%%%%%%%%%%%%%%%%%%%%%%%%%%%%%%%%%%%%%%%%%%%%%%%%%%%%%%%%%%%%%%%%

The discrete model is posed on a given finite set of vertices $\mathcal{V}$ and a set of unoriented edges (vessels)  $\mathcal{I}$.
We assume that the unoriented graph $(\mathcal{V},\mathcal{I})$ is connected and each pair of vertices $i$, $j\in\mathcal{V}$ is linked by at most one edge $(i,j)\in\mathcal{I}$.
Furthermore, $Q_{ij}, L_{ij}$ and, resp., $C_{ij}$ denote the flow through, length of and, resp., conductivity
of the edge $(i,j)\in \mathcal{I}$.
For biological applications, the Reynolds number of the flow is typically small and the flow is predominantly in the laminar (Poiseuille) regime.
Then the flow through the edge $(i,j)\in \mathcal{I}$ is proportional to the conductance $C_{ij}/L_{ij}$ and pressure drop along the edge, i.e.,
\begin{align}\label{flux}
   Q_{ij} = C_{ij} \frac{P_i-P_j}{L_{ij}}.
\end{align}
Conservation of mass in each vertex is expressed in terms of the Kirchhoff law
\begin{align}\label{kirchhoff}
    \sum_{j\in N(i)} Q_{ij} =  \sum_{j\in N(i)} C_{ij} \frac{P_i-P_j}{L_{ij}} = S_i \qquad \text{for all } i\in \mathcal{V},
\end{align}
where $N(i)$ denotes the set of edges adjacent to vertex $i\in\mathcal{V}$
and $S=(S_i)_{i\in\mathcal{V}}$ is the prescribed strength of the flow source ($S_i>0$) or sink ($S_i<0$)
at vertex $i\in\mathcal{V}$. Note that the flux is oriented, i.e. the flow rate from vertex $i$ to vertex $j$ is denoted by $Q_{ij}$ and it holds that $Q_{ij} =-Q_{ji}$.
We assume conservation of total mass, i.e.,
\begin{align*}
   \sum_{i\in \mathcal{V}} S_i = 0,
\end{align*}
which is a necessary (but not sufficient) condition for the solvability of \eqref{kirchhoff}.
For a given vector of conductivities $(C_{ij})_{(i,j)\in \mathcal{I}}$, the linear system \eqref{kirchhoff} for $(P_i)_{i\in \mathcal{V}}$
has a solution (unique up to an additive constant) if the underlying graph is connected,
where only edges with positive conductivities are taken into account
(i.e., edges with zero conductivities are discarded); see, e.g., \cite{gross2013handbook}.

Assuming that the material (metabolic) cost for an edge $(i,j)\in \mathcal{I}$ is proportional to a power $C_{ij}^\gamma$ of its conductivity $C_{ij}>0$, $\gamma \geq 0$,
Hu and Cai \cite{hu2013adaptation} considered the energy functional
\begin{align}\label{energy_discrete}
   E[C] = \sum_{\substack{(i,j)\in \mathcal{I}\\i< j}} \left( \frac{Q_{ij}[C]^2 }{C_{ij}} + \frac{\nu}{\gamma} C_{ij}^\gamma \right) L_{ij},
\end{align}
where $i<j$ in the summation symbol means that we count each edge only once.
The first term corresponds to the pumping power (kinetic energy) needed to pump the material through an edge $(i,j)\in \mathcal{I}$,
while the second term is the metabolic energy needed to maintain the edge, with $\nu>0$ the metabolic coefficient.
Hu and Cai considered the formal gradient flow of \eqref{energy_discrete} constrained by the Kirchhoff law \eqref{kirchhoff}, which is given by the ODE system
\begin{align}\label{micro_system}
   \frac{\mathrm{d}C_{ij}}{\mathrm{d}t} = \left(\frac{Q_{ij}[C]^2 }{C_{ij}^2}-\nu C_{ij}^{\gamma-1}\right) C_{ij} L_{ij},
\end{align}
see \cite{HKM} for details of the derivation.
%Note that the Kirchhoff law \eqref{kirchhoff} implies that $\frac{\partial Q_i(C)}{\partial C_i}=0$.
They provided a numerical evidence that the optimal networks generated by \eqref{micro_system},
i.e., minimizers of the energy \eqref{energy_discrete}, possess a phase transition at $\gamma = 1$:
for $\gamma>1$ the optimal network is tiled with loops, while for $\gamma <1$ it is a loopless tree.
In the next section we provide an analytical proof of the fact that for $\gamma <1$ the energy minimizer
does not contain any loops.
See \cite{Banavar} for a proof of an analogous statement for a model where the metabolic cost is given
as a constraint for the pumping energy. Let us note that the proof of \cite{Banavar} uses fundamentally different
techniques compared to our proof below, which is based on the concavity of a relaxed energy functional
(pressureless formulation, see also Section \ref{sec:pressureless}).

%%%%%%%%%%%%%%%%%%%%%%%%%%%%%%%%%%%%%%%%%%%%%%%%%%%%%%%%%%%%%%%%%%%%%%%%%%%%%%%%%%%%%%%%%%%
\subsection{No loops for \boldmath{$\gamma < 1$}}\label{subsec:noloops}
%%%%%%%%%%%%%%%%%%%%%%%%%%%%%%%%%%%%%%%%%%%%%%%%%%%%%%%%%%%%%%%%%%%%%%%%%%%%%%%%%%%%%%%%%%%
%In this section we prove that the energy minimizer of contains no loops, i.e., is a tree. For this sake, we consider a relaxed problem where
For a given $S=(S_i)_{i\in\mathcal{I}}$ we consider the relaxed problem
where the energy functional
\begin{align}\label{ECQ}
   \widetilde E[C, Q] = \sum_{\substack{(i,j)\in \mathcal{I}\\i< j}} \left( \frac{Q_{ij}^2 }{C_{ij}} + \frac{\nu}{\gamma} C_{ij}^\gamma \right) L_{ij}
\end{align}
is constrained by the local mass conservation law
\begin{align}  \label{localmass}
  \sum_{j\in N(j)} Q_{ij} = S_i \quad \mbox{ for all } j\in \mathcal{V}.
\end{align}
I.e., we consider the pressureless formulation where we drop the requirement
that fluxes are defined in terms of pressure differences as in \eqref{flux}.
The pressures $(P_j)_{j\in \mathcal{V}}$ can be recovered as Lagrange multipliers
in the constrained minimization problem \eqref{ECQ}--\eqref{localmass}.
Note that $Q_{ij}>0$ means net flow of material from the vertex $i$ to vertex $j$.
Since in \eqref{ECQ} the conductivities and fluxes are decoupled,
we can first minimize with respect to $C$, i.e., define
\begin{align}   \label{F}
   F[Q] := \inf_{C\in\R^N_+} \widetilde E[C,Q],
\end{align}
where $N\in\N$ denotes the number of edges. The minimizing vector of conductivities has the entries $C_{ij} = |Q_{ij}|^\frac{2}{\gamma+1}$ for all $(i,j)\in \mathcal{I} $.
Thus, defining
\begin{align} \label{f}
   f_\gamma(s) := (\gamma+1) |s|^\frac{2\gamma}{\gamma+1},
\end{align}
we have
\begin{align}  \label{F2}
   F[Q] = \sum_{\substack{(i,j)\in \mathcal{I}\\i< j}} f_\gamma(Q_{ij}) L_{ij}.
\end{align}
The existence of minimizers of $F=F[Q]$ on $\R^N$ follows
trivially from its continuity, coercivity and boundedness from below.
Moreover, note that for $0< \gamma < 1$ the function $f_\gamma$ given by \eqref{f} is strictly concave on the intervals $(-\infty,0)$ and $(0,\infty)$.
%while for $\gamma>1$ it is strictly convex on $\R$.
\vspace{2mm}

\begin{lemma}\label{lem:treeF}
Let $0 < \gamma < 1$.
Then the global minimizer $Q\in\R^N$ of the functional \eqref{F} constrained by the local mass conservation \eqref{localmass}
does not contain any loops, i.e., there exists no closed circle of edges $\widetilde{\mathcal{I}} := \{ (i_1,i_2), (i_2,i_3), \ldots, (i_K, i_1) \} \subset \mathcal{I}$
such that the fluxes $Q_{i_1 i_2}, Q_{i_2 i_3}, \dots, Q_{i_K i_1}$ are all nonzero.
\end{lemma}

\begproof
For contradiction, assume that there exists such a loop in the minimizer $Q$.
%, i.e., a closed chain of edges $\bar E \subset E$ with all nonzero fluxes $Q_i\neq 0$, $i\in\widetilde E$.
Moreover, we may assume without loss of generality that the loop $\widetilde{\mathcal{I}}$ is minimal,
i.e., the vertices belonging to $\widetilde{\mathcal{I}}$ are connected to each other exclusively
by the edges $(i_1,i_2), (i_2,i_3), \ldots, (i_K, i_1)$ and no other direct connection (edge) exists among them.
Then we have
\[
   \bar \eps := \min_{(i,j)\in\widetilde{\mathcal{I}}} |Q_{ij}| > 0.
\]
Moreover, let us assume, without loss of generality, that the vertices belonging to the loop
are indexed such that $i_1 < i_2 < \dots < i_K$.
Choosing any $\eps\neq 0$ such that $|\eps|<\bar\eps$, we define the modified fluxes
\begin{align}
   \widetilde Q_{ij} := \left\{ \begin{aligned}
      Q_{ij} + \eps\sign(i-j) \qquad\mbox{for } (i,j)\in\widetilde{\mathcal{I}},\\
      Q_{ij}\qquad\mbox{for } (i,j)\in \mathcal{I}\setminus \widetilde{\mathcal{I}}.
      \end{aligned} \right.
\end{align}
Then we have $\sign(\widetilde Q_{ij}) = \sign(Q_{ij})$ for all edges $(i,j)\in \mathcal{I}$,
i.e., we are introducing a (small) circular flux along $\widetilde{\mathcal{I}}$, but we are
not changing the sign (i.e., direction) of the individual fluxes $Q_{ij}$.

Let us show the local mass conservation \eqref{localmass} is still valid for the modified fluxes $\widetilde Q_{ij}$:
For a node $i$ belonging to the loop $\widetilde{\mathcal{I}} $ we have
%obviously, \eqref{localmass} is modified only for nodes adjacent to edges from $\widetilde E$.
%Since $\widetilde E$ is a loop, every such node $k\in\V$ has exactly two adjacent edges $j$, $\ell$ from $\widetilde E$,
\[
   \sum_{j\in N(i)} \widetilde Q_{ij} = \sum_{j\in N(i)\cap \widetilde{\mathcal{I}}} \widetilde Q_{ij} + \sum_{j\in N(i) \setminus \widetilde{\mathcal{I}}} Q_{ij}.
\]
Since $\widetilde{\mathcal{I}}$ is a minimal loop, there are exactly two edges, $(i,k)$ and $(i,\ell)$ in $\widetilde{\mathcal{I}}$
that are adjacent to the vertex $i$. Therefore,
\begin{align*}
   \sum_{j\in N(i)\cap \widetilde{\mathcal{I}}} \widetilde Q_{ij}  = \widetilde{Q}_{ik}+ \widetilde{Q}_{i\ell}
     &= Q_{ik} +  \eps\sign(i-k) + Q_{i\ell} + \eps\sign(i-\ell).
\end{align*}
Due to the above assumption about indexing, we have either $k < i < \ell$ or $\ell < i < k$.
In either case, $\sign(i-k) + \sign(i-\ell) = 0$ and, consequently,
\begin{align*}
   \sum_{j\in N(i)\cap \widetilde{\mathcal{I}}} \widetilde Q_{ij}  &= \widetilde{Q}_{ik}+ \widetilde{Q}_{i\ell} = Q_{ik} + Q_{i\ell},
\end{align*}
so that we indeed have the local mass conservation
\[
   \sum_{j\in N(i)} \widetilde Q_{ij} = \sum_{j\in N(i)} Q_{ij}  = S_i.
\]
Then, by \eqref{F2} and the definition of $\widetilde Q_{ij}$, we have
\[
   F[\widetilde Q] - F[Q] = \sum_{\substack{(i,j)\in\widetilde{\mathcal{I}}\\ i<j}} \bigl( f_\gamma(\widetilde Q_{ij}) - f_\gamma(Q_{ij}) \bigr) L_{ij}
        = \sum_{\substack{(i,j)\in\widetilde{\mathcal{I}}\\ i<j}} \bigl( f_\gamma(Q_{ij} - \eps) - f_\gamma(Q_{ij}) \bigr) L_{ij}.
\]
Since, by construction, $\sign(Q_{ij} - \eps) = \sign(Q_{ij})$,
and due to the strict concavity of $f_\gamma$ on $(-\infty,0)$ and $(0,\infty)$,
we have
\[
   f_\gamma(Q_{ij} - \eps) - f_\gamma(Q_{ij}) < - f_\gamma'(Q_{ij})\eps,
\]
so that, for a suitable choice of $\eps$ (positive or negative),
\[
   F[\widetilde Q] - F[Q] < - \eps \sum_{\substack{(i,j)\in\widetilde{\mathcal{I}}\\ i<j}} f_\gamma'(Q_{ij}) \leq 0.
\]
We conclude that if a loop of nonzero fluxes exists in $Q$,
then $Q$ cannot be the minimizer of $F$ subject to the local mass conservation \eqref{localmass}.
\endproof
\vspace{2mm}

The above Lemma states that global minimizers of $F=F[Q]$ given by \eqref{F} with $0 < \gamma < 1$,
subject to the local mass conservation \eqref{localmass}, are trees.
In fact, since the functional $F=F[Q]$ is strictly concave, and the constraint \eqref{localmass}
represents a polyhedral set, the minimizers are located in its extremal points.
Thus, structures with loops cannot be extremal points.
A direct consequence of Lemma \ref{lem:treeF} is that also global minimizers of $E=E[C]$ given by \eqref{energy_discrete},
constrained by the Kirchhoff law \eqref{kirchhoff}, are trees:
\vspace{2mm}

\begin{theorem}\label{thm:treeE}
Let $0 < \gamma < 1$.
Then the global minimizer $C\in\R^N_+$ of the functional $E = E[C]$ given by \eqref{energy_discrete} constrained by the Kirchhoff law \eqref{kirchhoff}
does not contain any loops, i.e., there exists no connected circle of edges $\widetilde{\mathcal{I}} := \{ (i_1,i_2), (i_2,i_3), \ldots, (i_K, i_1) \} \subset \mathcal{I}$
such that the fluxes $Q_{i_1 i_2}, Q_{i_2 i_3}, \dots, Q_{i_K i_1}$ are all nonzero.
\end{theorem}

\begproof
Obviously, for a fixed graph and a fixed distribution $S=(S_i)_{i\in\mathcal{I}}$ of sinks and sources, we have
\begin{align}
   \min_{C\in\R^N_+} \left\{ E[C] \mbox{ defined by \eqref{energy_discrete} with $Q=Q[C]$ subject to \eqref{kirchhoff}} \right\} \nonumber\\
   \geq
   \min_{C\in\R^N_+,\, Q\in\R^N} \left\{ \widetilde E[C,Q] \mbox{ defined by \eqref{ECQ} with $Q$ subject to \eqref{localmass}} \right\} \nonumber\\
   = \min_{Q\in\R^N} \left\{ F[Q] \mbox{ defined by \eqref{F} with $Q$ subject to \eqref{localmass}} \right\}.
\end{align}
If we prove that every minimizer $Q\in\R^N$ of $F=F[Q]$ is also a solution of the Kirchhoff law \eqref{kirchhoff},
i.e., construct pressures $P_i$ such that $Q_{ij} = C_{ij}\frac{P_i-P_j}{L_{ij}}$
with $C_{ij} = |Q_{ij}|^\frac{2}{\gamma+1}$, then this vector of conductivities %$C= (C_{ij})$
is a minimizer of the energy $E=E[C]$ defined by \eqref{energy_discrete}.
Then, by Lemma \ref{lem:treeF} we conclude that $C=(C_{ij})$ represents a tree.

The construction of the pressures is easy: Since $Q$, being a global minimizer of $F = F[Q]$,
is a tree (Lemma \ref{lem:treeF}), then we can start in any node $j_0$ and set $P_{j_0}:=0$.
Then we proceed inductively, for every edge $(j_0,j)$ emanating from $j_0$, we set
$C_{j_0 j} = |Q_{j_0 j}|^\frac{2}{\gamma+1}$ and
\[
   P_j := P_{j_0} - \frac{Q_{j_0 j} L_{j_0 j}}{C_{j_0 j}}.
\]
If some flux $Q_{j_0 j}$ happens to be zero, we set $C_{j_0 j}:=0$ and $P_j := P_{j_0}$.
Obviously, proceeding inductively in $j\in\mathcal{V}$, after a finite number of steps we assign conductivities to all edges
and pressures to all nodes in the graph,
verifying the Kirchhoff law \eqref{kirchhoff}.
\endproof
\vspace{2mm}

\begin{remark}
Suppose that $S$ consists only of one source and one sink, i.e., there are
two vertices $j_0$, $j_1\in \mathcal{V}$ such that
\[
   S_{j_0} = 1, \qquad S_{j_1} = -1,
\]
and $S_j=0$ for all other vertices $j\in \mathcal{V}$.
Then, by Theorem \ref{thm:treeE}, the global minimizer $C$
of the energy $E=E[C]$ with $0<\gamma<1$ is a path $\bar{\mathcal{J}} \subset \mathcal{J} $
connecting the vertices $j_0$ and $j_1$.
Then, obviously, all fluxes along the path are equal to $1$ or $-1$,
and the conductivities of edges belonging to $\bar{\mathcal{J}}$ are all equal to some $\bar C>0$.
Then
\[
   E[C] = \left( \frac{1}{\bar C} + \frac{\nu}{\gamma} \bar C^\gamma \right)  \sum_{\substack{(i,j)\in\bar {\mathcal{J}}\\ i<j}} L_{ij}.
\]
Consequently, the energy minimizer represents the shortest path connecting $j_0$ and $j_1$
in the sense of minimal sum of edge lengths $L_{ij}$.
\end{remark}
\vspace{2mm}

\begin{remark}
Note that the minimization problem for $F=F[Q]$ given by \eqref{F2} with $f(s) := (\gamma+1) |s|^\frac{2\gamma}{\gamma+1}$, i.e.,
\[
   \min_{Q\in\R^N} \sum_{(i,j)\in\mathcal{I}} |Q_{ij}|^\frac{2\gamma}{\gamma+1} L_{ij},
\]
constrained by \eqref{localmass}, coincides with the \emph{irrigation problem} introduced by Gilbert in \cite{Gilbert}.
There the energy is assumed to be a concave function of the flux (or amount of transported material), i.e.,
$0 < \frac{2\gamma}{\gamma+1} < 1$ which corresponds to $0 < \gamma<1$.
%For $\gamma=0$ we obtain the \emph{Steiner's minimal length problem}.. but we consider the edge lengths to be fixed.
The functional $F=F[Q]$ is then called the \emph{cost of irrigation} or \emph{Gilbert energy}.
In Section 4.4.2 of \cite{santambrogio2015optimal} it is mentioned (but not proved) that minimizers of the (concave) irrigation problem do not contain cycles.
\end{remark}

%%%%%%%%%%%%%%%%%%%%%%%%%%%%%%%%%%%%%%%%%%%%%%%%%%%%%%%%%%%%%%%%%%%%%%%%%%%%%%%%%%%%%%%%%%%
\subsection{Loops for \boldmath{$\gamma > 1$}}\label{subsec:notrees}
%%%%%%%%%%%%%%%%%%%%%%%%%%%%%%%%%%%%%%%%%%%%%%%%%%%%%%%%%%%%%%%%%%%%%%%%%%%%%%%%%%%%%%%%%%%
It is tempting to expect that the global minimizer of the energy functional \eqref{energy_discrete} with $\gamma>1$
will have "as many loops as possible". However, let us demonstrate with a simple example that the presence or absence of loops
in the energy minimizer depends on the data, in particular, the configuration of sources/sinks $S=(S_i)_{i\in\mathcal{V}}$ and the edge lengths.
%the transportation properties of the network, given by the configuration of sources/sinks $S=(S_i)_{i\in\mathcal{V}}$.
\vspace{2mm}

\begin{example}\label{ex:1}
Let us set $\gamma=2$ and consider a network consisting of three nodes $\mathcal{V} = \{1,2,3\}$ and two (unoriented) edges $\mathcal{I}=\{(1,2),(2,3)\}$ of unit length.
The sources and sinks be given by
\[
   S_1=-1,\qquad S_2=2,\qquad S_3=-1.
\]
Then the local mass conservation \eqref{localmass} gives
\[
   Q_{12} = 1,\qquad Q_{32} = 1,
\]
and the energy \eqref{F2} with $\gamma=2$ is equal to $F[Q] = 3(Q_{12}^2 + Q_{32}^2) = 6$.

Now, let us insert the third edge $(1,3)$ with a nonvanishing flux $q\neq 0$, closing the loop.
Then we have the new fluxes
\[
   \widetilde Q_{12} := 1+q, \qquad \widetilde Q_{32} = 1-q,\qquad \widetilde Q_{31} = q,
\]
obviously satisfying the local mass conservation \eqref{localmass},
and the corresponding energy is
\[
   F[\widetilde Q] = 3\left( (1+q)^2 + (1-q)^2 + q^2 \right) = 6 + 3q^2 > F[Q].
\]
We conclude that, even though $\gamma>1$, closing the loop by edge insertion with nonvanishing flux leads to energy increase.
\end{example}
\vspace{2mm}

\begin{example}\label{ex:2}
We keep the setting of the previous example, i.e., $\gamma=2$, $\mathcal{V} = \{1,2,3\}$, $\mathcal{I}=\{(1,2),(2,3)\}$,
but now the sources and sinks be given by
\[
   S_1=-1,\qquad S_2=3,\qquad S_3=-2.
\]
The fluxes imposed by the local mass conservation \eqref{localmass} read then
\[
   Q_{12} = 1,\qquad Q_{32} = 2,
\]
and the corresponding energy $F[Q] = 3(Q_{12}^2 + Q_{13}^2) = 15$.

Now, let us again insert the third edge $(1,3)$ with the flux $q\neq 0$.
The new fluxes are
\[
   \widetilde Q_{12} := 1+q, \qquad \widetilde Q_{32} = 2-q,\qquad \widetilde Q_{31} = q,
\]
obviously satisfying the local mass conservation \eqref{localmass},
and the corresponding energy is
\[
   F[\widetilde Q] = 3\left( (2-q)^2 + (1+q)^2 + q^2 \right) = 3(5 - 2q + 3q^2),
\]
which is strictly less than $F[Q] = 15$ if $0<q<2/3$ (and attaining minimum for $q=1/3$).
Thus, in this case the creation of a loop does reduce the energy consumption of the network.
\end{example}

Note that the setting in Example \ref{ex:1} was symmetric, i.e., same edge lengths
and both sinks of the same intensity. In this case, the optimal structure is a tree even for $\gamma>1$.
In Example \ref{ex:2} we disturbed the symmetry by considering sinks of different intensity.
Alternatively, we could have consider different edge lengths, which would also lead
to the conclusion that closing the loop reduces the energy consumption.

%%%%%%%%%%%%%%%%%%%%%%%%%%%%%%%%%%%%%%%%%%%%%%%%%%%%%%%%%%%%%%%%%%%%%%%%%%%%%%%%%%%%%%%%%%%%
\section{Mesoscopic modeling approach for transportation networks}\label{sec:mesoscopic}
%%%%%%%%%%%%%%%%%%%%%%%%%%%%%%%%%%%%%%%%%%%%%%%%%%%%%%%%%%%%%%%%%%%%%%%%%%%%%%%%%%%%%%%%%%%%
In this section we shall discuss the mesoscopic modeling approach, briefly introduced in \cite{albi2017continuum}.
The key idea is to interpret the term
\[
   \frac{Q_{ij}^2 }{C_{ij}^2} = \left(\frac{P_j-P_i}{L_{ij}}\right)^2
\]
in the microscopic system \eqref{micro_system} as a difference quotient for a continuous pressure variable
$p=p(x)$ in the direction $\theta_{ij}\in \mathbb{S}^1_+$ of the edge $(i,j)\in\mathcal{I}$, i.e.,
\[
   \frac{P_j-P_i}{L_{ij}} \approx \theta_{ij} \cdot \nabla p(x_{ij}),
\]
where $x_{ij}$ denotes the midpoint of the edge $(i,j)$ and $\mathbb{S}^1_+$ is the unit half-sphere in $\R^d$, i.e.,
\begin{align}
   \mathbb{S}_+^1 = \{\theta \in \R^d;\; |\theta|=1, \theta_1\geq 0\}.
\end{align}
Note that we may restrict $\theta$ to the half-sphere since the edges are not oriented.
Then, we introduce the continuum conductivity $C=C(x,\theta)$ for $x\in\Omega$, with $\Omega$ a bounded domain in $\R^d$
where the discrete graph is embedded, and $\theta\in \mathbb{S}^1_+$.
Identifying $C_{ij}=C_{ij}(t)$ with $C(t,x_{ij}, \theta_{ij})$, equation \eqref{micro_system} can be rewritten as
\begin{align}\label{equ5}
   \frac{dC}{dt}(x_{ij},\theta_{ij}) =  \bigl( c_0^2 |\theta_{ij}\cdot \nabla p(x_{ij})|^2 - C(x_{ij},\theta_{ij})^{\gamma-1} \bigr)  C(x_{ij},\theta_{ij}),
\end{align}
where we introduced $c_0^2>0$, the \emph{activation constant}, and set $\nu=1$.
Interpreting \eqref{equ5} as a characteristic system for a hyperbolic transport equation for
a probability measure $\mu_t=\mu_t(x,\theta,C)$, we obtain
\begin{align}\label{meso_system1}
   \partial_t \mu_t + \partial_C \left( (c_0^2 |\theta \cdot \nabla p|^2 -C^{\gamma-1}) C \mu_t \right) = 0,
\end{align}
where $\mu_t=\mu_t(x,\theta,C)$ represents the probability to have an edge of conductivity $C\geq 0$ in direction $\theta\in\mathbb{S}^1_+$
at $x\in\Omega$ and time $t\geq 0$.
%in $\Omega \times \mathbb{S}_+^1\times \mathbb{R}_+$
The continuum version of the Kirchhoff law \eqref{kirchhoff} is then given by the Poisson equation
\begin{align}\label{meso_system2}
   -\nabla \cdot (\mathbb{P}[\mu]\nabla p)=S
\end{align}
with the permeability tensor   
\begin{align}\label{meso_system3}
   \mathbb{P}[\mu] = \int_\mathbb{R_+}\int_{\mathbb{S}_+^1} C \, \theta \otimes \theta \, \d\mu(\cdot \, , \theta, C),
\end{align}
where $S=S(x)$ is a given density of sources/sinks verifying the global mass conservation
\begin{align*}
   \int_\Omega S(x) \d x = 0.
\end{align*}
Obviously, the global mass conservation is a necessary condition for solvability of \eqref{meso_system2}
subject to the no-flux boundary condition
\begin{align*}
   n(x) \cdot \mathbb{P}[\mu]\nabla p(x) = 0\quad\mbox{for } x\in\partial\Omega,
\end{align*}
which we adopt in the sequel; $n(x)$ denotes the unit normal vector in $x\in\partial\Omega$.

Let us point out that, in general, the permeability tensor \eqref{meso_system3} is a matrix-valued measure on $\Omega$.
Consequently, well posedness of weak solutions of the Poisson equation \eqref{meso_system2} represents an open (and difficult) problem.
%in particular, \eqref{meso_system2} with a measure valued $\mathbb{P}[\mu]$ requires a proper reinterpretation.
In Section \ref{subsec:special} we present several special cases where we are able to provide a proper reinterpretation of \eqref{meso_system2}
such that at least its formal solvability is established. Development of a general well posedness theory for \eqref{meso_system2}--\eqref{meso_system3}
is beyond the scope of this paper.

Let us also remark that the system \eqref{meso_system1}--\eqref{meso_system3} has the formal Wasserstein-type gradient flow structure
\begin{align} \label{meso_system_gf}
   \partial_t \mu_t + \partial_C(C\mu_t(-\partial_C \mathcal{E}'[\mu_t]))=0,
\end{align}
with the energy functional
\begin{align} \label{Emu}
   \mathcal{E}[\mu] = \int_\mathbb{R_+}\int_{\mathbb{S}_+^1}\int_\Omega\left( c_0^2 C|\theta \cdot \nabla p[\mu]|^2 + \frac{C^\gamma}{\gamma} \right)\, \d \mu(x,\theta,C), 
\end{align}
where $p[\mu]$ solves \eqref{meso_system2} with no-flux boundary conditions and given source density $S=S(x)$ independent of time.
The first term in the energy functional corresponds to the network-fluid interaction energy,
while the second one is the metabolic (relaxation) energy.

%%%%%%%%%%%%%%%%%%%%%%%%%%%%%%%%%%%%%%%%%%%%%%%%%%%%%%%%%%%%%%%%%%%%%
\subsection{Connection to the macroscopic model of \cite{hu2013optimization}}
%%%%%%%%%%%%%%%%%%%%%%%%%%%%%%%%%%%%%%%%%%%%%%%%%%%%%%%%%%%%%%%%%%%%%
Based on heuristic modeling arguments, Hu and Cai proposed in \cite{hu2013optimization} a macroscopic network formation model of the form
\begin{align}
   -\nabla \cdot [(m\otimes m)\nabla p] &= S,\label{macro_system1}\\
   \frac{\partial m}{\partial t} - c_0^2 (m\cdot \nabla p)\nabla p + |m|^{2(\gamma-1)}m &= 0,\label{macro_system2}
\end{align}
where $m=m(t,x)\in \mathbb{R}^d$ is the vector-valued local conductance of the network structure.
To establish a link between the mesoscopic model \eqref{meso_system1}--\eqref{meso_system2}
and the macroscopic model \eqref{macro_system1}--\eqref{macro_system2} we construct the conductance vector $m$ as
\begin{align*}
   m(t,x) = \int_{\mathbb{R}_+}\int_{\mathbb{S}_+^1} \sqrt{C}\theta \mu_t(x,  \d\theta,\d C),
\end{align*}
and adopt the monokinetic closure
\begin{align*}
   \mu_t(x,\theta,C) = \delta(C-\hat{C}(t,x)) \otimes \delta(\theta-\hat{\theta}(t,x))
\end{align*}
for suitable functions $\hat{C}=\hat{C}(t,x)$ and $\hat{\theta}=\hat{\theta}(t,x)$.
We obtain that $m(t,x)=\sqrt{\hat{C}(t,x)}\hat{\theta}(t,x)$ and from \eqref{meso_system2} that
\[
    \mathbb{P}[\mu] =  \hat{C} \hat{\theta}\otimes\hat{\theta} = m\otimes m.
\]
Consequently, \eqref{macro_system1} is satisfied. 
Moreover, using the identities $\hat{C}=|m|^2$ and $\hat{\theta}=\frac{m}{|m|}$ and formal integration by parts, we obtain
\begin{align*}
   \frac{\partial m}{\partial t} &= \frac{\partial}{\partial t} \int_{\mathbb{R}_+}\int_{\mathbb{S}_+^1} \sqrt{C}\theta \mu_t(\cdot,  \d\theta, \d C)\\
   &=- \int_{\mathbb{R}_+}\int_{\mathbb{S}_+^1} \sqrt{C}\theta \partial_C \left[ (c_0^2 |\theta\cdot \nabla p|^2 -C^{\gamma-1}) C \mu_t(\cdot,  \d\theta, \d C) \right]\\
   &=\frac{1}{2} \int_{\mathbb{R}_+}\int_{\mathbb{S}_+^1} \sqrt{C}\theta \, (c_0^2 |\theta\cdot \nabla p|^2 -C^{\gamma-1}) \mu_t(\cdot,  \d\theta, \d C)\\
   &=\frac{1}{2} \left(c_0^2\left|\frac{m}{|m|}\cdot \nabla p\right|^2 - |m|^{2(\gamma-1)}\right) m.
\end{align*}
Hence, we obtain an equation similar to \eqref{macro_system2}. However, while the activation term in \eqref{macro_system2}
is proportional to $\nabla p$, here it is a scalar multiple of $m$. This is due to the fact that in the microscopic model the locations
and directions of the edges are fixed (and given a priori). Thus, the direction of an edge never changes;
only its conductivity is adapted according to the transportation needs of the network.
The local direction of the edge is on the PDE level described by the vector $m=m(t,x)$.
The fact that the right-hand side above is a scalar multiple of $m$ means that the
direction of $m$ does not change during the evolution - in agreement with the discrete model.
Only the length of $m$ evolves in time, which corresponds to adaptation of the conductivity of the corresponding edge.
%Consequently, there is a mismatch between the heuristic PDE system of \cite{hu2013optimization} and the formal macroscopic closure
%of the mesoscopic model \eqref{meso_system1}--\eqref{meso_system2}, which in turn directly reflects the microscopic model of \cite{hu2013adaptation}
%described in section \ref{sec:discrete}.

Let us note that a regularized version of the PDE system \eqref{macro_system1}--\eqref{macro_system2} of the form
\begin{align*}
   -\nabla \cdot [(rI+m\otimes m )\nabla p]&=S, \\
   \frac{\partial m}{\partial t}-D^2 \Delta m -c_0^2 (m\cdot \nabla p)\nabla p +|m|^{2(\gamma-1)}m&=0,
\end{align*}
with the constant diffusivity $D\geq 0$ and the isotropic background permeability $r = r(x)\geq r_0 > 0$,
has been studied in the series of papers \cite{haskovec2015mathematical, haskovec2016notes, albi2016biological, albi2017continuum}.

%%%%%%%%%%%%%%%%%%%%%%%%%%%%%%%%%%%%%%%%%%%%%%%%%%%%%%%%%%%%%%%%%%%%%%%%%%%%%%%%%%%%%%%%%%%%
\subsection{Special stationary solutions.}\label{subsec:special}
%%%%%%%%%%%%%%%%%%%%%%%%%%%%%%%%%%%%%%%%%%%%%%%%%%%%%%%%%%%%%%%%%%%%%%%%%%%%%%%%%%%%%%%%%%%%
In this Section, we discuss some special stationary solutions of system \eqref{meso_system1}--\eqref{meso_system3}.
Besides providing a connection to another well-studied problem by Putti et al, cf. \cite{facca2018towards},
we show that solutions of the discrete model \eqref{micro_system} are special stationary solutions of the mesoscopic one.
However, as we aim to construct solutions concentrated on a lower dimensional subdomain,
we need a proper reinterpretation of the Poisson equation \eqref{meso_system2} as the following discussion will show. 

%%%%%%%%%%%%%%%%%%%%%%%%%%%%%%%%%%%%%%%%%%%%%%%%%%%%%%%%%%%%%%%%%%%%%%%%%%%%%%%%%%%%%%%%%%%%
\subsubsection{Edges aligned with the pressure gradient - system of Putti et al.}\label{sssec:Putti}
%%%%%%%%%%%%%%%%%%%%%%%%%%%%%%%%%%%%%%%%%%%%%%%%%%%%%%%%%%%%%%%%%%%%%%%%%%%%%%%%%%%%%%%%%%%%
Based on the intuition that energetically optimal networks should have edges aligned with the pressure gradients of the transported material,
we introduce the ansatz
\begin{align}  \label{ansatz1}
   \mu(x,\theta,C) = \lambda(x,C) \otimes \delta\left(\theta - {\frac{\nabla p}{|\nabla p|}(x)}\right),
\end{align}
where $\lambda=\lambda(x,C)$ is a probability measure on $\Omega\times\R^+$ and
$\delta$ denotes the Dirac delta.
Inserting \eqref{ansatz1} into \eqref{meso_system3} gives
\[
   \mathbb{P}[\mu] = \bar{C} \frac{\nabla p \otimes \nabla p}{|\nabla p|^2}   \qquad\mbox{with } \bar{C}:= \int_0^\infty C \d\lambda(\cdot,C),
\]
and the Poisson equation \eqref{meso_system2} reads
\[
   - \nabla \cdot (\bar{C}(x) \nabla p) = S.
\]
The stationarity condition for \eqref{meso_system1} reads
\[
   c_0^2 C|\theta\cdot\nabla p|^2 - C^\gamma = 0 \qquad\mbox{on supp}(\mu)
\]
and with \eqref{ansatz1},
\[
   c_0^2 C|\nabla p|^2 - C^\gamma = 0 \qquad\mbox{on supp}(\lambda).
\]
Making the further ansatz $\lambda(x,C) = \delta(x-\bar C(x))$ we obtain
\begin{align*}
    c_0^2 \bar{C}|\nabla p|^2 - \bar{C}^\gamma &= 0 \\
    -\nabla \cdot (\bar{C} \nabla p)&=S,
\end{align*}
which for $\gamma=1$ reduces to the system proposed and studied by Putti et al. \cite{facca2018towards}, %, system (1.2) in \cite{facca2018towards},
which with our notation reads
\begin{align}\label{Putti}
\begin{aligned}
    (c_0^2 |\nabla p|^2 - 1) \bar{C} &= 0 \\
    -\nabla \cdot (\bar{C} \nabla p)&=S.
\end{aligned}
\end{align}
%the problem is analogous to Putti et al, cf. \cite{facca2018towards}, where they study existence and uniqueness.
This system is proposed in \cite{facca2018towards} as a formal extension to the continuous setting of a discrete model \cite{Tero} describing
the dynamics of slime mold, Physarum Polycephalum, which was used to simulate the
ability of the slime mold to find the shortest path connecting two food sources in a maze. The discrete model
describes the dynamics of the slime mold on a finite-dimensional planar graph. Similarly as the model of
Hu and Cai \cite{hu2013adaptation} it is using a pipe-flow analogy whereby mass transfer occurs because of pressure differences
with a conductivity coefficient that varies with the flow intensity.
The authors of \cite{facca2018towards} proved the well-posedness of the system \eqref{Putti} for small times.
Moreover, they provide motivations for the conjecture that the system presents a time-asymptotic equilibrium which
is a solution of Monge-Kantorovich partial differential equations governing optimal transportation problems;
see, e.g., \cite{evans1999differential, ambrosio2003lecture}.

%%%%%%%%%%%%%%%%%%%%%%%%%%%%%%%%%%%%%%%%%%%%%%%%%%%%%%%%%%%%%%%%%%%%%%%%%%%%%%%%%%%%%%%%%%%%
\subsubsection{Single source and single sink.}\label{ssec:singlesourcesink}
%%%%%%%%%%%%%%%%%%%%%%%%%%%%%%%%%%%%%%%%%%%%%%%%%%%%%%%%%%%%%%%%%%%%%%%%%%%%%%%%%%%%%%%%%%%%
In this section we focus on the situation of having a single point source $x^+\in\Omega$ and a single point sink in $x^-\in\Omega$ described by
the signed measure
\[
   S(x) =\delta(x-{x^+}) - \delta(x-{x^-}).
\]
Let $\Gamma: (0,1) \to \Omega$ be any open smooth curve connecting $x^+$ and $x^-$
and $\mathcal{t}: \Gamma \to \mathbb{S}^1$ its tangent vector.
%By a slight abuse of notation, we denote $\mathcal{t}(x)$ the tangent vector at $x\in\Gamma$.
Then we claim that any measure $\mu$ of the form
\begin{align} \label{ansatz2}
   \mu(x,\theta,C) &= \delta_\Gamma(x) \otimes \delta(\theta-\mathcal{t}(x)) \otimes \delta(C-\bar C) +\eta (x, \theta ) \otimes \delta(C)
\end{align}
with $\bar{C}=c_0^{\frac{2}{1+\gamma}}$ and $\eta$ an arbitrary probability measure on $\Omega\times\mathbb{S}^1$
is a stationary solution of \eqref{meso_system1}--\eqref{meso_system3}.
Here $\delta_\Gamma$ denotes the one-dimensional Hausdorff measure concentrated on the curve $\Gamma$.

Indeed, plugging the ansatz \eqref{ansatz2} in the definition of the permeability tensor given in \eqref{meso_system2}, we obtain
\begin{align}\label{perm_measure}
\mathbb{P}[\mu](x)= \bar{C} \, \mathcal{t}(x)\otimes \mathcal{t}(x) \, \delta_\Gamma(x).
\end{align}
%Note that the following calculations are just formal. As the permeability tensor \eqref{perm_measure} is measure valued, the analysis of the Poisson equation is not covered by the standard elliptic theory.
Let $s\in(0,1)$ be the arclength variable along the curve $\Gamma$.
A formal integration by parts with a test function $\phi\in C^\infty(\Omega)$,
assuming continuity of $\nabla p$ on $\Omega$, yields
\begin{align*}
   \int_\Omega \nabla \cdot (\mathbb{P}[\mu] \nabla p)\phi \, \d x
      &= -\int_\Omega \bar{C} (\mathcal{t}\otimes \mathcal{t} \nabla p)\cdot \nabla \phi \, \delta_\Gamma(x)\,\d x \\
      &=-\int_\Omega \bar{C} (\mathcal{t} \cdot \nabla p) (\mathcal{t} \cdot \nabla \phi) \,\delta_\Gamma(x)\, \d x\\
      &=-\int_0^1 \bar{C} (\partial_s p) (\partial_s \phi) \, \d s\\
      &= \int_0^1 \partial_s (\bar{C}\partial_s p) \phi \, \d s.\\
\end{align*}
Thus, only the values of $\bar{C}$ and $\nabla p$ along $\Gamma$ are relevant, and we have the formal identity
\begin{align}\label{equ2}
\nabla \cdot (\mathbb{P}[\mu] \nabla p) &= \partial_s (\bar{C}\partial _s p)
\end{align}
along the curve $\Gamma$.
As the sources/sinks are only concentrated in the nodes $x^+$ and $x^-$,
we have that $S=0$ in the interior of the curve $\Gamma$ and, consequently,
\begin{align*}
   \partial_{s}(\bar{C} \partial_{s} p) \equiv 0.
\end{align*}
Thus, $\bar{C} \partial_{s} p\equiv G$ for some constant $G$ on $\Gamma$. This leads to
\begin{align*}
   \int_\Omega \mathbb{P}[\mu] \nabla p \cdot \nabla \phi \,\delta_\Gamma(x)\, \d x
   &= \int_0^1 \bar{C} \partial_s p \, \partial_s \phi \,\d s\\
   &= G \int_0^1 \partial_s \phi(x(s))\, \d s = G(\phi(x^-)-\phi(x^+)).
\end{align*}
On the other hand, since
\begin{align*}
   \int_\Omega S\phi\, \d x = \phi(x^+) - \phi(x^-),
\end{align*}
we have $G=1$ and $\bar C \partial_s p \equiv 1$.
With \eqref{ansatz2} the stationarity condition for \eqref{meso_system1} reads
\begin{align}   \label{stat2}
   c_0^2 C|\mathcal{t}\cdot\nabla p|^2 - C^\gamma = 0 \qquad\mbox{on }\Gamma,
\end{align}
so that
\begin{align*}
   \bar{C}^\gamma &= c_0^2 \bar{C} | \mathcal{t}\cdot \nabla p|^2 = c_0^2 \bar{C} |\partial_s p|^2 = c_0^2/\bar C. % \quad \text{for } x\in \Gamma.
\end{align*}
Consequently, \eqref{stat2} is satisfied with $\bar{C}=c_0^{\frac{2}{1+\gamma}}$.

Observe that with the constant conductivity $\bar{C}$ and constant derivative $\partial_s p \equiv 1$ along $\Gamma$
the energy \eqref{Emu} only depends on the length of the curve $\Gamma$.
Consequently, the global energy minimizer realizes the shortest path connecting
the source $x^+\in\Omega$ and the sink $x^-\in\Omega$.
In this context, let us again refer to the experiments manifesting the ability of slime mold to find shortest connections in a maze \cite{facca2018towards, Tero}.

%%%%%%%%%%%%%%%%%%%%%%%%%%%%%%%%%%%%%%%%%%%%%%%%%%%%%%%%%%%%%%%%%%%%%%%%%%%%%%%%%%%%%%%%%%%%
\subsubsection{Discrete network solution.}\label{ssec:discretenetsol}
%%%%%%%%%%%%%%%%%%%%%%%%%%%%%%%%%%%%%%%%%%%%%%%%%%%%%%%%%%%%%%%%%%%%%%%%%%%%%%%%%%%%%%%%%%%%
Generalizing the example from the previous section, we show that a discrete network structure realized by a connected graph
with a finite set of edges and vertices can be interpreted as a stationary solution to the mesoscopic model \eqref{meso_system1}.
Let the graph consist of a finite set $\mathcal{V}$ of nodes connected by a finite set $\mathcal{I}$ of edges as described in Section \ref{sec:discrete}.
In particular, each edge $i\in \mathcal{I}$ is an open, smooth nonintersecting curve $\Gamma_i$ in $\R^d$, also not intersecting any other edge.
For $x\in\Gamma_i$ we denote $\mathcal{t}_i(x)\in\mathbb{S}^1$ the tangent vector to $\Gamma_i$ at $x$.
We then define the set $\mathcal{G}\subset \Omega\times\mathbb{S}^1\times \R_+$,
\[
   \mathcal{G} := \left\{ (x,\theta,C);\; x\in \Gamma_i \mbox{ for some } i\in\mathcal{I},\, \theta=\mathcal{t}_i(x),\, C=C_i \right\},
\]
where $C_i>0$ denotes the conductivity of the edge $i\in \mathcal{I}$.
Let $\mu = \mu(x,\theta,C)$ be the Hausdorff measure on $\Omega\times\mathbb{S}^1\times \R_+$ concentrated on the set $\mathcal{G}$.
We claim that with a suitable choice of the conductivities $C_i$, $\mu$ can be formally interpreted as a stationary solution to \eqref{meso_system1}--\eqref{meso_system3}.
As in Section \ref{ssec:singlesourcesink}, we need to reinterpret the Poisson equation involving the measure valued permeability tensor
\begin{align*}
   \mathbb{P}[\mu](x) = \int_{\R_+}\int_{\mathbb{S}_+^1} C \, \theta \otimes \theta \, \mu(x, \d\theta, \d C)
      = \sum_{i\in\mathcal{I}} {C}_i \, \mathcal{t}_i(x)\otimes\mathcal{t}_i(x) \, \delta_{\Gamma_i}(x),
\end{align*}
where $\delta_{\Gamma_i}$ is the one-dimensional Hausdorff measure in $\R^d$ concentrated on the curve $\Gamma_i$.

The weak formulation of the left-hand side of \eqref{meso_system1} with a test function $\phi\in C^\infty(\Omega)$ becomes
\begin{align*}
   \int_\Omega \mathbb{P}[\mu]\nabla p\cdot \nabla \phi\, \d x
      &=\sum_{i\in \mathcal{I}} \int_{\Gamma_i} {C}_i (\mathcal{t}_i(x) \cdot \nabla p) (\mathcal{t}_i(x) \cdot \nabla \phi)\,\d \delta_{\Gamma_i}(x)\\
      &=\sum_{i\in \mathcal{I}} \int_0^1 (\partial_{s_i} p) (\partial_{s_i} \phi) \, \d s_i,
\end{align*}
where $s_i\in(0,1)$ is the arc length parameter of the curve $\Gamma_i$.
The sources/sinks are concentrated on the nodes, in particular,
\[
   S(x) = \sum_{j\in\mathcal{V}} S_j \delta(x-x_j),
\]
where $x_j\in\R^d$ is the location of the node $j\in\mathcal{V}$ and we impose the global mass conservation
\[
   \sum_{j\in\mathcal{V}} S_j = 0.
\]
Consequently, we have $\partial_{s_i}({C}_i \partial_{s_i} p) \equiv 0$ in the interior of each edge $\Gamma_i$,
and therefore  ${C}_i \partial_{s_i} p\equiv G_i$ for some constants $G_i$, $i\in \mathcal{I}$.
Hence,
\begin{align*}
   \int_\Omega \mathbb{P}[\mu]\nabla p\cdot \nabla \phi\, \d x
      &= \sum_{i\in \mathcal{I}} G_i\int_0^1 \partial_{s_i} \phi \,\d s_i\\
      &= \sum_{i\in \mathcal{I}}G_i(\phi(x_i^-)-\phi(x_i^+))
      = \sum_{j\in \mathcal{V}} S_j\phi(x_j),
\end{align*}
where $x_i^+$ and $x_i^-$ denote the endpoints of the curve $\Gamma_i$. 
Choosing $\phi$ supported on a small neighborhood of the node $j\in\mathcal{V}$, we obtain
\begin{align*}
   \sum_{i\in N(j)}G_i\phi(x_j) &= S_j\phi(x_j),
\end{align*}
where $N(j)$ denotes the set of edges adjacent to the node $j$.
Consequently,
\begin{align}  \label{sumG}
   \quad \sum_{i\in N(j)} G_i = S_j.
\end{align}
Moreover, the fact that $\partial_{s_i} p$ is constant on each edge implies
\begin{align*}
   \partial_{s_i} p = \frac{p(x_i^-)-p(x_i^+)}{L_i}\qquad\mbox{for } s_i\in (0,1),
\end{align*}
where $L_i$ is the length of the edge $i\in\mathcal{I}$.
Using the identity ${C}_i \partial_{s_i} p\equiv G_i$ in \eqref{sumG} finally gives
\begin{align*}
   \sum_{i\in N(j)} {C}_i \frac{p(x_i^-)-p(x_i^+)}{L_i}&= S_j \text{ for all } j\in \mathcal{V},
\end{align*}
which corresponds to the Kirchhoff law \eqref{kirchhoff}.

The stationarity condition for \eqref{meso_system1} reads
\begin{align}   \label{stat3}
   c_0^2 C_i |\mathcal{t}_i\cdot\nabla p|^2 - C_i^\gamma = 0 \qquad\mbox{on }\Gamma_i,
\end{align}
so that
\begin{align*}
   {C}_i^{\gamma+1} = c_0^2 {C}_i^2 |\partial_{s_i}p|^2 = c_0^2 G_i^2 \qquad\mbox{for } i\in\mathcal{I}.
\end{align*}
Consequently, \eqref{stat3} is verified by choosing ${C}_i=(c_0 G_i)^{\frac{2}{\gamma +1}}$,
where $(G_i)_{i\in\mathcal{I}}$ solves \eqref{sumG}.
This shows the consistency of our mesoscopic model \eqref{meso_system1} with the microscopic model
of Hu and Cai \cite{hu2013adaptation} described in Section \ref{sec:discrete}.

%%%%%%%%%%%%%%%%%%%%%%%%%%%%%%%%%%%%%%%%%%%%%%%%%%%%%%%%%%%%%%%%%%%%%%%%%%%%%%%%%%%%%%%%%%%%%%%%%%%%%
\section{Monokinetic model for the conductivity}\label{sec:monokinetic}
%%%%%%%%%%%%%%%%%%%%%%%%%%%%%%%%%%%%%%%%%%%%%%%%%%%%%%%%%%%%%%%%%%%%%%%%%%%%%%%%%%%%%%%%%%%%%%%%%%%%%
Based on the expectation that equation \eqref{meso_system1} tends to produce concentrations in the conductivity variable,
we propose the ansatz
\begin{align}\label{mono_capacity}
   \mu_t(x,\theta,C) = \eta(x,\theta) \delta(C-\hat{C}(t,x,\theta)),
\end{align}
where $\eta=\eta(x,\theta)$ is a probability measure on $\Omega\times\mathbb{S}^1_+$.
Note that $\eta$ does not depend on time. This reflects the fact that the discrete model of Section \ref{sec:discrete}
is set on an a-priori given graph (network structure), with fixed locations of vertices and directions of edges.
This network structure is encoded by the time-independent measure $\eta=\eta(x,\theta)$.
Taking the first moment of \eqref{meso_system1} with respect to $C$,
\[
   \frac{\d}{\d t} \int_\mathbb{R_+} C \mu(\,\cdot\, , \, \cdot\, , \d C)
      = \int_\mathbb{R_+}(c_0^2 C|\theta \cdot \nabla p|^2 -C^\gamma )\, \mu(\,\cdot\, , \, \cdot\, , \d C),
\]
and inserting the ansatz \eqref{mono_capacity}, we obtain
\begin{align}\label{meso2}
   \partial_t \hat{C} = c_0^2 \hat{C} |\theta\cdot \nabla p|^2 -\hat{C}^\gamma\qquad\mbox{on supp}(\eta).
\end{align}
This is coupled to the Poisson equation
\begin{align}\label{meso2_1}
   -\nabla \cdot (\mathbb{P}[\hat{C}] \nabla p) = S,
\end{align}
subject to the no-flux boundary condition on $\partial\Omega$,
and with the permeability tensor
\begin{align}\label{meso2_2}
   \mathbb{P}[\hat{C}] = \int_{\mathbb{S}_+^1}(\theta \otimes \theta) \hat{C} \, \d\theta.
\end{align}
For notational convenience we shall omit the hat in $\hat{C}$ in the sequel.
Note that the system \eqref{meso2}-\eqref{meso2_1} exhibits the formal gradient flow structure
with respect to the Fisher-Rao metric, cf. \cite{gallouet2017jko,chizat2016interpolating} as well as Section \ref{sec:pressureless}, i.e.
\[
   \partial_t C = - C \partial_C \mathcal{E}[C],
\] 
with the energy functional
\begin{align} \label{Emono}
   \mathcal{E}[C] = \int_{\mathbb{S}_+^1}\int_\Omega c_0^2 C|\theta \cdot \nabla p|^2 +\frac{C^\gamma}{\gamma} \d x \, \d\theta,
\end{align}
where $p=p(x)$ solves \eqref{meso2_1}.

Below we shall construct particular stationary solutions of the system \eqref{meso2}--\eqref{meso2_1}.

%%%%%%%%%%%%%%%%%%%%%%%%%%%%%%%%%%%%%%%%%%%%%%%%%%%%%%%%%%%%%%%%%%%%%%%%%%%%%%%%%%%%%%%%%%%%%%%%%%%%%
\subsection{Stationary solutions of \eqref{meso2}--\eqref{meso2_1} in one spatial dimension}
%%%%%%%%%%%%%%%%%%%%%%%%%%%%%%%%%%%%%%%%%%%%%%%%%%%%%%%%%%%%%%%%%%%%%%%%%%%%%%%%%%%%%%%%%%%%%%%%%%%%%
Let us consider the spatially one-dimensional setting $\Omega=(0,1)$.
Then the permeability tensor \eqref{meso2_2} reduces to $\mathbb{P}[C]=C$.
Assuming $C \geq C_0 >0$ on $\Omega$ and taking into account the no-flux boundary condition for the Poisson equation \eqref{meso2_1} at $x=0$, i.e., $p'(0)=0$,
we integrate \eqref{meso2_1} to obtain
\begin{align*}
   p'(x)&= \frac{B(x)}{C(x)} \qquad\mbox{with } B(x) = \int_0^x S(y)\, \d y.
\end{align*}
Thus, \eqref{meso2} reads
\begin{align}\label{stat1D}
   \partial_t C&=\frac{c_0^2 B^2}{C}-C^\gamma,
\end{align}
and we obtain the stationary state 
\begin{align}\label{stat1D_2}
   C_\infty(x)&=(c_0 B(x))^{\frac{2}{\gamma+1}}.
\end{align}
Denoting the right side of the ODE \eqref{stat1D} by $f(C):=\frac{c_0^2 B^2}{C}-C^\gamma$, differentiating and plugging in the stationary state \eqref{stat1D_2} yields
\begin{align*}
   f'(C_\infty(x))&=-(1+\gamma)(c_0^2B(x)^2)^{\frac{\gamma-1}{\gamma+1}} \leq 0.
\end{align*}
Therefore, $C_\infty(x)$ is asymptotically stable as long as $B(x)\neq 0$.
This result can be compared with Section 6.1 of \cite{haskovec2015mathematical}, where the stationary solutions of the spatially one-dimensional
version of the model proposed by \cite{hu2013optimization} have been studied. This model exhibits a larger variety of stationary solutions with stability
depending on the value of the parameter $\gamma>0$, with a bifurcation at $\gamma=1$.

%%%%%%%%%%%%%%%%%%%%%%%%%%%%%%%%%%%%%%%%%%%%%%%%%%%%%%%%%%%%%%%%%%%%%%%%%%%%%%%%%%%%%%%%%%%%%%%%%%%%%
\subsection{Stationary solutions of \eqref{meso2}-\eqref{meso2_1} for \boldmath{$\gamma>1$}}\label{ss1}
%%%%%%%%%%%%%%%%%%%%%%%%%%%%%%%%%%%%%%%%%%%%%%%%%%%%%%%%%%%%%%%%%%%%%%%%%%%%%%%%%%%%%%%%%%%%%%%%%%%%%
Inserting the nonzero stationary solution
$C=\left(c_0^2|\theta \cdot \nabla p|^2 \right)^{\frac{1}{\gamma-1}}$
of \eqref{meso2} into the Poisson equation \eqref{meso2_1} and taking into account the no-flux boundary condition,
we obtain the weak formulation
\begin{align}\label{stat_poisson}
    \int_\Omega \int_{\mathbb{S}_+^1} \left(c_0^2|\theta \cdot \nabla p|^2 \right)^{\frac{1}{\gamma-1}}(\theta\cdot\nabla p)(\theta\cdot\nabla\phi) \d\,\theta \d x = \int_\Omega S \phi \d x
\end{align}
for all test functions $\phi\in C^\infty(\Omega)$.
Solutions of this nonlinear elliptic equation for $\gamma>1$ can be constructed using the direct method of calculus of variations:

\begin{proposition}%{(Existence of network patterns)\\}
Let $\gamma>1$.
For every $S\in L^2(\Omega)$ there exists a unique solution $p\in H^1(\Omega)$ of \eqref{stat_poisson}
satisfying $\int_\Omega p(x) \d x = 0$.
\end{proposition}

\begproof
We define the functional $\mathcal{F}:H^1(\Omega)\to \mathbb{R}$ by
\begin{align}\label{el_stat_poisson}
   \mathcal{F}[p]:=\int_\Omega\int_{\mathbb{S}_+^1}\frac{\gamma-1}{2\gamma-1} c_0^{\frac{2}{\gamma-1}}|\theta\cdot\nabla p|^{\frac{2\gamma-1}{\gamma-1}}\, d\theta \, dx-\int_\Omega pS\, dx
\end{align}
and $\mathcal{F}[p]:=\infty$ if $\nabla p\notin L^{\frac{2\gamma-1}{\gamma-1}}$.
Since $\frac{2\gamma-1}{\gamma-1}>2$, the functional $\mathcal{F}$ is uniformly convex on $H^1(\Omega)$.
Moreover, a straightforward application of the Poincar\'{e} inequality provides coercivity on the set $H^{1,0}(\Omega):=\left\{ p\in H^1(\Omega);\; \int_\Omega p(x) \d x = 0 \right\}$.
The classical theory (see, e.g. \cite{evans1997partial}) provides the existence of a unique minimizer $p\in H^{1,0}(\Omega)$.
Since, as is easy to check, \eqref{el_stat_poisson} is the Euler-Lagrange equation for the minimizer of $\mathcal{F}$,
$p$ is the unique solution in $H^{1,0}(\Omega)$ of \eqref{stat_poisson}.
\endproof

%%%%%%%%%%%%%%%%%%%%%%%%%%%%%%%%%%%%%%%%%%%%%%%%%%%%%%%%%%%%%%%%%%%%%%%%%%%%%%%%%%%%%%%%%%%%%%%%%%%%%
\subsection{Stationary solutions of \eqref{meso2}-\eqref{meso2_1} with \boldmath{$\gamma=1$}}\label{ss2}
%%%%%%%%%%%%%%%%%%%%%%%%%%%%%%%%%%%%%%%%%%%%%%%%%%%%%%%%%%%%%%%%%%%%%%%%%%%%%%%%%%%%%%%%%%%%%%%%%%%
The stationary solution of \eqref{meso2} with $\gamma=1$ satisfies
\begin{align}\label{meso3}
    C(c_0^2|\theta\cdot \nabla p|^2 -1) = 0,
\end{align}
i.e.,
\begin{align}\label{meso4}
   c_0|\theta\cdot \nabla p(x)| =1\qquad \mbox{on supp}(C(x,\cdot)) \mbox{ for all } x\in\Omega.
\end{align}
Clearly, \eqref{meso4} is solvable for $\theta\in\mathbb{S}^1_+$ if and only if $c_0 |\nabla p(x)|\geq 1$ on supp$(C(x,\cdot))$.
To avoid the complexities of the construction of a general solution, we shall propose, under an additional structural assumption on the function $S=S(x)$,
the particular solution
\begin{align}\label{stat_sol_2d}
   C(x,\theta) = \alpha(x)\delta\left(\theta - {\hat{\theta}}(x)\right),
\end{align}
with suitable $\alpha=\alpha(x)$ and $\hat{\theta} = \hat{\theta}(x)$ such that
$|c_0\nabla p| \equiv 1$ on $\Omega$.
Let us assume that $S=S(x)$ is of the form
\[
   S = - \nabla\cdot w\quad\mbox{in }\Omega,\qquad w\cdot\nu = 0\quad\mbox{on }\partial\Omega,
\]
with an irrotational vector field $w=w(x)$.
Inserting the ansatz \eqref{stat_sol_2d} into the Poisson equation \eqref{meso2_1} gives then
\begin{align}  \label{CrazyPoisson}
    \left(\chi^+(x) - \chi^-(x)\right) \frac{1}{c_0} \alpha(x) \hat{\theta}(x) = w(x),
\end{align}
where $\chi^\pm = \chi^\pm(x)$ are the characteristic functions of the disjoint sets $A^+$, $A^-$
such that $A^+ \cup A^- = \mbox{supp}(w)$.
We then set
\begin{align}  \label{meso5}
   \nabla p:= \frac{1}{c_0} \frac{w}{|w|},\qquad
   \hat\theta := \left(\chi^+ - \chi^-\right) \frac{w}{|w|},\qquad
   \alpha := c_0 |w|
\end{align}
on $\mbox{supp}(w)$. Note that $\nabla\times w = 0$ implies $\nabla\times\frac{w}{|w|} = 0$,
so that $\frac{w}{|w|}$ is indeed a gradient field.
For a given $w=w(x)$ we choose the sets $A^+$, $A^-$
such that $\hat\theta=\hat\theta(x)$ defined above is an element of the hemisphere $\mathbb{S}^1_+$ for all $x\in\Omega$.
It is easily checked that \eqref{stat_sol_2d} with \eqref{meso5} is a solution of \eqref{meso2}--\eqref{meso2_1}.

Note that the ansatz \eqref{stat_sol_2d} for $C=C(x,\theta)$ can be interpreted
as a network of infinitesimal tubes (edges) with local direction of the tube given by $\hat{\theta}(x)$
and local conductivity $\alpha(x)$.
Observe that similarly to Section \ref{sssec:Putti} the directions of the tubes (edges) $\hat{\theta}(x)$
are aligned with the pressure gradients.

%%%%%%%%%%%%%%%%%%%%%%%%%%%%%%%%%%%%%%%%%%%%%%%%%%%%%%%%%%%%%%%%%%%%%%%%%%%%%%%%%%%%%%%%%%%%%%%%%%%%%
\section{Pressureless formulation of the monokinetic model}\label{sec:pressureless}
%%%%%%%%%%%%%%%%%%%%%%%%%%%%%%%%%%%%%%%%%%%%%%%%%%%%%%%%%%%%%%%%%%%%%%%%%%%%%%%%%%%%%%%%%%%%%%%%%%%%%
In this section we shall introduce a relaxed formulation of the system \eqref{meso2}--\eqref{meso2_2} that avoids the explicit presence of the pressure variable.
In this way we can overcome the solvability issues in the Poisson equation \eqref{meso2_2} due to the possible degeneracy of the permeability tensor \eqref{meso2_1}.
We introduce the scalar variable $Q=Q(t,x,\theta)$ such that $Q\theta$ is the flux in direction $\theta\in\mathbb{S}^1_+$, i.e., 
we identify $Q=C \theta \cdot \nabla p$.
We then consider the relaxed version of the energy functional \eqref{Emono},
\begin{align}\label{energy_pless}
   \widetilde{\mathcal{E}}[C,Q] := \int_{\mathbb{S}_+^1} \int_\Omega c_0^2 \frac{Q^2}{C}+\frac{C^\gamma}{\gamma} \d x \d\theta,
\end{align}
subject to the linear constraint
\begin{align}\label{linC}
   -\nabla \cdot \left(\int_{\mathbb{S}_+^1} Q \theta \d\theta \right) = S\qquad\mbox{on } \Omega.
\end{align}
Imposing the no-flux boundary condition on $\partial\Omega$, \eqref{linC} is written in the weak form as
\begin{align}\label{linearC}
   \int_\Omega \int_{\mathbb{S}_+^1} \theta\cdot\nabla\phi \d Q(x, \theta) = \int_\Omega S \phi \d x\qquad
   \mbox{for all } \phi\in C^\infty(\Omega),
\end{align}
allowing for measure-valued $Q$.
Note that the relation $Q=C \theta \cdot \nabla p$ is obtained as the formal optimality condition
for the constrained minimization problem
\begin{align} \label{Q[C]}
    \min_{Q} \widetilde{\mathcal{E}}[C,Q]\quad \text{ subject to } \eqref{linearC},
   %-\nabla \cdot \left(\int_{\mathbb{S}_+^1} Q \theta \d\theta \right) = S\qquad\mbox{on } \Omega,
\end{align}
where the pressure $p=p(x)$ is recovered as a multiple of the Lagrange multiplier.

We shall consider the formal constrained gradient flow of \eqref{energy_pless}--\eqref{linearC}
with respect to the Fisher-Rao metric, cf. \cite{gallouet2017jko,chizat2016interpolating},
\begin{align}\label{flux_eq1}
   \partial_t C = \frac{c_0^2 Q[C]^2}{C}-C^\gamma,
\end{align}
where $Q[C]$ is the minimizer in \eqref{Q[C]}. Note that \eqref{Q[C]} is a convex minimization problem
with the convex constraint \eqref{linearC}, so that a unique $Q[C]$ exists for each $C\geq 0$.

The Fisher-Rao distance between nonnegative measures $C_0, C_1\in \mathcal{M}_+(\Omega\times\mathbb{S}^1_+)$ is defined as
\begin{align}\label{dist1}
   \mathcal{d}_\mathrm{FR}(C_0,C_1)^2 &:= \inf \left\{ \int_0^1 \int_\Omega\int_{\mathbb{S}^1_+} |r_t(x)|^2 \,\d C_t(x,\theta) \,\d t; \;
     {(C,r)\in \mathcal{A}[C_0,C_1]} \right\},
\end{align}
where $\mathcal{A}[C_0,C_1]$ consists of curves $[0,1]\ni t\mapsto (C_t,r_t)\in \mathcal{M}_+(\Omega\times \mathbb{S}_+^1)\times L^2(\Omega\times \mathbb{S}_+^1, \d C_t)$
such that $t\mapsto C_t$ is narrowly continuous with endpoints $C_0$, $C_1$ and
\[
   \partial_t C_t =C_t r_t
\]
in the sense of distributions on $(0,1)\times \Omega\times \mathbb{S}_+^1$.
Moreover, if $\lambda\in \mathcal{M}_+(\Omega\times\mathbb{S}^1_+)$ is any reference measure such that $C_0$, $C_1$ are absolutely continuous with respect to $\lambda$,
with Radon-Nikodym derivatives $\frac{\d C_0}{\d\lambda}$, $\frac{\d C_1}{\d\lambda}$, then we have
\begin{align}\label{dist2}
   \mathcal{d}_\mathrm{FR}(C_0,C_1)^2  &= 4 \int_\Omega \int_{\mathbb{S}_+^1} \left|\sqrt{\frac{\d C_0}{\d\lambda}}-\sqrt{\frac{\d C_1}{\d\lambda}}\right|^2\, \d\lambda(x,\theta).
\end{align}
Note that by $1$-homogeneity, the above expression does not depend on the choice of $\lambda$.
%Hence, the distance measures an energy of growth $\|\tilde{r}\|_{L^2(dC)}$. In the case of smooth densities, we have the relation $\tilde{r}=\frac{r}{C}$.
The gradient with respect to the Fisher-Rao metric of a functional $\mathcal{F}: \mathcal{M}_+(\Omega\times\mathbb{S}^1_+) \to \R$ at $C\in\mathcal{M}_+(\Omega\times\mathbb{S}^1_+)$
is $C \mathcal{F}'(C)$, where $\mathcal{F}'(C)$ is the Fr\'{e}chet derivative, c.f. \cite{gallouet2017jko}. This explains why \eqref{flux_eq1} is the formal Fisher-Rao
gradient flow of \eqref{energy_pless}--\eqref{linearC}.

%%%%%%%%%%%%%%%%%%%%%%%%%%%%%%%%%%%%%%%%%%%%%%%%%%%%%%%%%%%%%%%%%%%%
\subsection{Reformulation of the energy functional \eqref{energy_pless} with $\gamma=1$}\label{ssec:ReformE}
For $\gamma=1$ it is possible to reformulate \eqref{energy_pless} such that its domain of definition is extended to the space
of measures. Obviously, we have to give a meaning to the term $Q^2/C$, which is possible by using the following
result related to the Benamou-Brenier functional of optimal transport theory \cite{santambrogio2015optimal}.
\vspace{2mm}

\begin{lemma}[Lemma 5.17 from \cite{santambrogio2015optimal}] \label{lem:BB}
Let $p, q>1$ be a pair of exponents such that $\frac{1}{p}+\frac{1}{q}=1$.
Set $K_q:=\{(a,b)\in \mathbb{R}^2:\, a+\frac{1}{q}|b|^q\leq 0\}$. Then for $(x,y)\in\mathbb{R}^2$ we have
\begin{align*}
\sup_{(a,b)\in K_q} \{ax+by\}=f_p(x,y):=
\begin{cases}
\frac{1}{p}\frac{|y|^p}{x^{p-1}}\quad &\text{if }x>0\\
0\quad&\text{if }x=0,\, y=0\\
+\infty\quad &\text{if } x=0,\, y\neq 0,\text{ or } x<0.
\end{cases}
\end{align*}
Moreover, being a supremum of linear functions, $f_p$ is convex and lower semicontinuous.
\end{lemma}

\noindent
Using the above Lemma with $p=q=2$, we can rewrite \eqref{energy_pless} with $\gamma=1$ as
\begin{align*}
   \widetilde{\mathcal{E}}[C,Q] = \int_{\mathbb{S}_+^1} \int_\Omega 2 c_0^2 f_2(C, Q) + C \d x \d\theta.
\end{align*}
Then, Proposition 5.18 of \cite{santambrogio2015optimal} states that if both $C$ and $Q$ are absolutely continuous with respect to
the Lebesgue measure on $\Omega\times \mathbb{S}_+^1$, we have
\begin{align*}
   \widetilde{\mathcal{E}}[C,Q] = \sup \left\{ \int_{\mathbb{S}_+^1} \int_\Omega 2 c_0^2 (aC + bQ) + C \d x \d\theta;\; (a,b)\in C_b(\Omega\times\mathbb{S}_+^1, K_2) \right\}.
\end{align*}
Thus, the definition of \eqref{energy_pless} can be extended to $C, Q\in \mathcal{M}(\Omega\times\mathbb{S}_+^1)$ as
\begin{align} \label{enM}
   \widetilde{\mathcal{E}}[C,Q] &= \sup \left\{ \int_{\mathbb{S}_+^1} \int_\Omega (2 c_0^2 a + 1) \d C(x,\theta)\right.\nonumber\\ &\qquad \qquad \qquad \qquad + \int_{\mathbb{S}_+^1} \int_\Omega b \d Q(x,\theta);\;
     \left. (a,b)\in C_b(\Omega\times\mathbb{S}_+^1, K_2) \right\}.
\end{align}
Note that the convexity and lower semicontinuity of the function $f_p$ asserted by Lemma \ref{lem:BB}
implies convexity and weak lower semicontinuity of $\widetilde{\mathcal{E}}[C,Q]$.

In the next section we shall make a first step towards developing a well posedness theory in measures
for the system \eqref{linearC}--\eqref{flux_eq1}, based on a minimizing movement scheme for the energy functional \eqref{enM}.

%%%%%%%%%%%%%%%%%%%%%%%%%%%%%%%%%%%%%%%%%%%%%%%%%%%%%%%%%%%%%%%%%%%%%%%%%%%%%%%%%%%%%%%%%%%%%%%%%%%%%
\subsection{A minimizing movement scheme for \eqref{linearC}--\eqref{flux_eq1}}\label{sec:min_mov}
%%%%%%%%%%%%%%%%%%%%%%%%%%%%%%%%%%%%%%%%%%%%%%%%%%%%%%%%%%%%%%%%%%%%%%%%%%%%%%%%%%%%%%%%%%%%%%%%%%%%%
Let us fix the time step $\tau>0$ and construct iteratively the sequence of pairs of measures $(C^n, Q^n)_{n\in\N}$
in $\mathcal{M}(\Omega\times\mathbb{S}_+^1)\times \mathcal{M}(\Omega\times\mathbb{S}_+^1)$,
\begin{align}\label{min_mov}
   (C^{n+1}, Q^{n+1}) := \argmin \Biggl\{ \frac{1}{2\tau}\mathcal{d}_\mathrm{FR} (C,C^n) + \widetilde{\mathcal{E}}[C,Q];\;
        C\in\mathcal{M}(\Omega\times\mathbb{S}_+^1),  \qquad \\
         Q\in\mathcal{M}(\Omega\times\mathbb{S}_+^1) \mbox{ verifying \eqref{linearC}} \Biggr\},  \nonumber
\end{align}
where $\widetilde{\mathcal{E}}[C,Q]$ is defined in \eqref{enM}.
The iterative scheme is subject to the initial datum $C^0\in \mathcal{M}_+(\Omega\times\mathbb{S}_+^1)$, $Q^0\in \mathcal{M}(\Omega\times\mathbb{S}_+^1)$
with $Q^0$ absolutely continuous with respect to $C^0$, $Q^0$ satisfying \eqref{linearC}, and having finite energy $\widetilde{\mathcal{E}}[C^0,Q^0] < +\infty$.
\vspace{2mm}

\begin{lemma}\label{lem:JKO}
The iterative scheme \eqref{min_mov} subject to the initial datum given above
admits a minimizing sequence $(C^n, Q^n)_{n\in\N}$ such that
\begin{align}\label{aprioriE}
   \widetilde{\mathcal{E}}[C^{n},Q^{n}] \leq \widetilde{\mathcal{E}}[C^0,Q^0] < +\infty \qquad\mbox{for all } n\in\N,
\end{align}
$C^n$ is a nonnegative measure on $\Omega\times\mathbb{S}_+^1$
for all $n\in\N$ and $Q^n$ is absolutely continuous with respect to $C^n$ for all $n\in\N$.
Moreover, the sequence $(C^n, Q^n)_{n\in\N}$ is uniformly bounded in
$\mathcal{M}_+(\Omega\times\mathbb{S}_+^1)\times \mathcal{M}(\Omega\times\mathbb{S}_+^1)$.
\end{lemma}

\begproof
The existence of the minimizer in \eqref{min_mov} for each $n\in\N$ is guaranteed by the convexity and weak lower semicontinuity
of the Fisher-Rao distance $\mathcal{d}_\mathrm{FR}$ and of the energy functional $\widetilde{\mathcal{E}}$.

By definition, the minimizing movement scheme \eqref{min_mov} yields
\[
   \frac{1}{2\tau} \mathcal{d}_\mathrm{FR}(C^{n+1},C^n)^2+{\cal E}(C^{n+1},Q^{n+1}) \leq {\cal E}(C^n,Q^n) \qquad\mbox{for all } n\in\N.
\]
Summing up with respect to $n$ directly implies \eqref{aprioriE}.

From Proposition 5.18 of \cite{santambrogio2015optimal} it follows that $\widetilde{\mathcal{E}}[C^{n},Q^{n}] < +\infty$
implies $C^n\geq 0$ and $Q^n$ absolutely continuous with respect to $C^n$.

Finally, denote $s=s(x,\theta)$ the sign of $Q^n$, i.e., $|s|\equiv 1$ on $\Omega\times\mathbb{S}_+^1$ and
\[
   \int_{\mathbb{S}_+^1}\int_\Omega s \d Q^n = \int_{\mathbb{S}_+^1}\int_\Omega \d |Q^n|,
\]
where $|Q^n|$ is the total variation of the measure $Q^n$.
Choose $b:=\frac{s}{\sqrt{2}c_0}$ and $a:= -\frac{1}{4c_0^2}$.
Then it is easy to check that $(a,b) \in K_2$ almost everywhere on $\Omega\times\mathbb{S}_+^1$,
with $K_2$ defined in Lemma \ref{lem:BB}.
By the nonnegativity of $C^n$ and \eqref{enM} we have then
\begin{align*}
   \widetilde{\mathcal{E}}[C^n,Q^n] \geq \frac12 \int_{\mathbb{S}_+^1} \int_\Omega \d |C^n|(x,\theta) + \frac{1}{\sqrt{2}c_0} \int_{\mathbb{S}_+^1} \int_\Omega \d |Q^n|(x,\theta).
\end{align*}
With \eqref{aprioriE} this implies the uniform boundedness of $(C^n, Q^n)_{n\in\N}$ in the space of measures.
\endproof

The uniform boundedness of the sequence $(C^n, Q^n)_{n\in\N}$ in the space of measures provided by Lemma \ref{lem:JKO}
implies weak* compactness for a subsequence of piecewise linear interpolates.
This facilitates the limit passage $\tau\to 0$ in the subsequence, yielding a limiting curve $(C_t, Q_t)_{t\in[0,T]}$
in the space of measures. However, due to the lack of regularity estimates, it is currently our of reach to prove that the limiting
curve is a measure valued solution of the system \eqref{linearC}--\eqref{flux_eq1}.

%%%%%%%%%%%%%%%%%%%%%%%%%%%%%%%%%%%%%%%%%%%%%%%%%%%%%%%%%%%%%%%%%%%%%%%%%%%%%%%%%%%%%%%%%%%%
\subsubsection{Single source and single sink.}\label{ssec:singlesourcesink2}
%%%%%%%%%%%%%%%%%%%%%%%%%%%%%%%%%%%%%%%%%%%%%%%%%%%%%%%%%%%%%%%%%%%%%%%%%%%%%%%%%%%%%%%%%%%%
In analogy to Section \ref{ssec:singlesourcesink} we construct a stationary solution of \eqref{linearC}--\eqref{flux_eq1}
for the situation when there is a single point source $x^+\in\Omega$ and a single point sink in $x^-\in\Omega$, i.e.,
\[
   S(x) =\delta(x-{x^+}) - \delta(x-{x^-}).
\]
Let $\Gamma: (0,1) \to \Omega$ be any open smooth curve connecting $x^+$ and $x^-$
and $\mathcal{t}: \Gamma \to \mathbb{S}^1$ its tangent vector.
In order to the construction be valid for general values of $\gamma>0$, we define the reference measure
\[
   \eta(x,\theta) := \delta_\Gamma(x) \otimes \delta(\theta-\mathcal{t}(x)),
\]
where $\delta_\Gamma$ denotes the one-dimensional Hausdorff measure concentrated on the curve $\Gamma$.
%and we again denoted $\mathcal{t}(x)$ the tangent vector at $x\in\Gamma$.
We then construct the measures $C$ and $Q$ to be absolutely continuous with respect to $\eta$,
with the densities (Radon-Nikodym derivatives)
\[
   \overline{C} := \frac{\d C}{\d\eta},\qquad \overline{Q} := \frac{\d Q}{\d\eta}.
\]
The stationary version of \eqref{flux_eq1} is then interpreted as an equation for the densities $\overline{C}$, $\overline{Q}$, i.e.,
\begin{align}\label{flux_eq1_bar}
   \frac{c_0^2 \overline{Q}^2}{\overline{C}}-\overline{C}^\gamma = 0\qquad\mbox{on supp}(\eta).
\end{align}
%We claim that the constant densities $\overline{C}=c_0^{\frac{2}{1+\gamma}}$ and $\overline{Q}=1$
%give a stationary solution of \eqref{flux_eq1_bar}, constrained by \eqref{linearC}.
Plugging in the ansatz for $Q$ into \eqref{linearC}, we have for all test functions $\phi\in C^\infty(\Omega)$,
\begin{align*}
   \int_\Omega \int_{\mathbb{S}_+^1} \theta\cdot\nabla\phi \d Q(x,\theta) = \int_\Omega \int_{\mathbb{S}_+^1} \theta\cdot\nabla\phi\; \overline{Q} \d \eta(x,\theta)
    = \int_\Gamma \mathcal{t}\cdot\nabla\phi\; \overline{Q} \d\Gamma(x).
\end{align*}
Changing to the arclength variable $s\in(0,1)$ along the curve $\Gamma$ and integrating formally by parts, we have
\begin{align*}
   \int_\Omega \int_{\mathbb{S}_+^1} \theta\cdot\nabla\phi \d Q(x,\theta) =
      \int_0^1 (\partial_s\phi) \overline{Q} \d s 
   & =- \int_0^1 (\partial_s\overline{Q}) \phi \d s \\
    & = \int_\Omega S \phi\, \d x.
\end{align*}
Since $S\equiv 0$ in the interior of the curve $\Gamma$, we have $\partial_s\overline{Q}\equiv 0$,
and similarly as in Section \ref{ssec:singlesourcesink} we conclude that $\overline{Q}\equiv 1$.
Inserting this into \eqref{flux_eq1_bar} implies finally that $\overline{C}=c_0^{\frac{2}{1+\gamma}}$.
Consequently, the energy \eqref{energy_pless} depends only on the length of the curve $\Gamma$,
so that the global energy minimizer realizes the shortest path connecting
the source $x^+\in\Omega$ and the sink $x^-\in\Omega$.

Let us note that the above construction can be extended to a system of curves (edges) connecting a finite number
of vertices where the sources/sinks are concentrated. The construction of a discrete network solution
would then follow the same steps as in Section \ref{ssec:discretenetsol}.

%%%%%%%%%%%%%%%%%%%%%%%%%%%%%%%%%%%%%%%%%%%%%%%%%%%%%%%%%%%%%%%%%%%%%%%%%%%%%%%%%%%%%%%%%%%%%%%%%%%%%
\subsection{Equivalence of energy minimization and Beckmann problem for \boldmath{$\gamma=1$}}
%%%%%%%%%%%%%%%%%%%%%%%%%%%%%%%%%%%%%%%%%%%%%%%%%%%%%%%%%%%%%%%%%%%%%%%%%%%%%%%%%%%%%%%%%%%%%%%%%%%%%
Finally, we again establish the connection between the minimizers of the energy functional \eqref{energy_pless} with $\gamma=1$
and the system \eqref{Putti} proposed in \cite{facca2018towards} to simulate the
ability of Physarum to find the shortest path connecting two food sources in a maze. The latter are solutions of the Beckmann problem, a well-known problem in optimal transport,
cf. \cite{santambrogio2015optimal}.

As in our reformulation below, we can write the energy for the Beckmann problem as
 \begin{align} \label{BeckmannF}
   \widetilde{\mathcal{F}}[\overline{C},q] = \sup \left\{  \int_\Omega (2 c_0^2 a + 1) \d \overline C(x) %\qquad \qquad \qquad \qquad
   +  \int_\Omega b \d  q(x);\;
      (a,b)\in C_b(\Omega, K_2) \right\}.
\end{align}
The Beckmann problem is then given by the minimization problem
\begin{equation} \label{Beckmann}
     \widetilde{\mathcal{F}}[\overline{C},q] \rightarrow \min \quad \text{ subject to } \quad - \nabla \cdot q = S.
\end{equation} 
We mention that for a minimizer the vector valued measure $q$ is absolutely continuous with respect to $\overline{C}$. 

On the one hand, any solution of the system \eqref{Putti} is an admissible minimizer of \eqref{energy_pless} subject to \eqref{linC}.

\begin{theorem}Let $(\overline{C}_*,q_*)$ be the  solution of the Beckmann problem \eqref{Beckmann} with $\frac{dq}{d\overline{C}} \neq 0$ $\overline{C}$-almost everywhere. Then $(C_*,Q_*)$ is a solution of \eqref{energy_pless} subject to \eqref{linC} if and only if
\begin{equation}
C_* = \overline{C}_* \delta_{\theta_*}, \qquad Q_* = q_* \cdot \theta_* \delta_{\theta_*}, \qquad \theta_* = \frac{q_*}{|q_*|}
\end{equation} 
hold $\overline{C}_*$-almost everywhere in $\Omega$. 
\end{theorem}
\begproof
Since $(\overline{C}_* \delta_{\theta_*},q_* \cdot \theta_* \delta_{\theta_*})$ is admissible for the minimization of $\widetilde{\mathcal{E}}$ we see that
$$ \widetilde{\mathcal{E}}[C_*,Q_*] \leq \widetilde{\mathcal{E}}[\overline{C}_* \delta_{\theta_*},q_* \cdot \theta_* \delta_{\theta_*}] = \widetilde{\mathcal{F}}[\overline{C}_*,q_*]. $$
On the other hand, given $C_*$ and $Q_*$ we construct admissible elements
$$ \overline{C} = \frac{1}{|\mathbb{S}_+^1|}  \int_{\mathbb{S}_+^1}  ~dC_*(\cdot,\theta), \qquad
q = \frac{1}{|\mathbb{S}_+^1|} \int_{\mathbb{S}_+^1} \theta ~dQ_*(\cdot,\theta), $$ 
and find by  convexity (the new variables arise as convex combinations of $C_*$ and $Q_*$) with $\hat{\theta} = \frac{q}{|q|}$
$$ \widetilde{\mathcal{F}}[\overline{C} ,q ] = \widetilde{\mathcal{E}}[\overline{C}  \delta_{\hat{\theta}},q \cdot \theta \delta_{\hat{\theta}}] \leq \widetilde{\mathcal{E}}[C_*,Q_*] . $$
Since $ \widetilde{\mathcal{F}}[\overline{C} ,q ]  \geq \widetilde{\mathcal{F}}[\overline{C}_* ,q_* ]$, we obtain the assertion. 
\endproof

\subsection*{Acknowledgements}
%%%%%%%%%%%%%%%%%%%%%%%%%%%%%%%%%%%%%%%%%%%%%%%%%%%%%%%%%%%%%%%%%%%%%%%%%%%%%%%%%%%%%
H.R. acknowledges support by the Austrian Science Fund (FWF) project F 65.
M.B. acknowledges support by ERC via Grant EU FP 7 - ERC Consolidator
Grant 615216 LifeInverse. M.B. would like to thank the Isaac Newton Institute for
Mathematical Sciences, Cambridge, for support and hospitality during the programme Variational
Methods for Imaging and Vision, where work on this paper was undertaken, supported
by EPSRC grant no EP/K032208/1 and the Simons foundation.
H.R. and J.H. would like to thank C. Schmeiser (U. Vienna) for several insightful discussions.

\bibliographystyle{plain}
\bibliography{bibliography}
\end{document}